\documentclass[10pt,twoside]{article}
\usepackage{amsfonts}
\usepackage{fancyhdr}
\usepackage{titlesec}
\usepackage{cite}
\usepackage{ifthen}
\usepackage{pifont}
\usepackage{stmaryrd}
\usepackage{setspace}
\usepackage{indentfirst}
\usepackage{amsmath,amssymb,amscd,bbm,amsthm,mathrsfs,dsfont}
\input amssym.def

\titleformat{\section}{\centering\large\bfseries}{\S\arabic{section}}{1em}{}
\newboolean{first}
\setboolean{first}{true}

\newtheorem{theorem}{Theorem}[section]
\newtheorem{lemma}{Lemma}[section]
\newtheorem{rem}{Remark}[section]

\newtheorem{definition}{Definition}[section]

\begin{document}

\setlength\abovedisplayskip{2pt}
\setlength\abovedisplayshortskip{0pt}
\setlength\belowdisplayskip{2pt}
\setlength\belowdisplayshortskip{0pt}

\title{\bf \Large Commutators of weighted Hardy operator on weighted $\lambda$-central Morrey space\author{Huihui Zhang$^1$ \ \ \ \  Yan Lin$^1$  \ \ \ \ Xiao Yu$^{2}$}\date{}} \maketitle

 \footnote{MR Subject Classification: 42B20, 42B25.}

 \footnote{Supported by the National Natural Science Foundation of China\,(11961056,11561057), the Jiangxi Natural Science Foundation of China\,(20192BAB2101004), the Science Foundation of Jiangxi Education Department\, (GJJ190890). }
\footnote{1 Department of Mathematics,  China University of Mining Technology(Beijing), Beiing 100083, P.R.China.}
\footnote{2 Department of Mathematics, Shangrao Normal University, Shangrao 334001, P.R.China.}
\footnote{Email address: zhanghuihuinb@163.com (H.Zhang);linyan@cumtb.edu.cn(Y.Lin); yx2000s@163.com(X.Yu).\\}
\begin{center}
\begin{minipage}{135mm}

{\bf \small Abstract}.\hskip 2mm {\small
In this paper, the authors prove the boundedness of commutators generated by the weighted Hardy operator on weighted $\lambda$-central Morrey space with the weight $\omega$ satisfying the doubling condition. Moreover,  the authors  give the characterization  for the  weighted $\lambda$-central Campanato space by introducing a new kind of operator which is related to the commutator of weighted Hardy operator.}
\end{minipage}
\end{center}

\thispagestyle{fancyplain} \fancyhead{}
\fancyhead[L]{\textit{}\\
} \fancyfoot{} \vskip 10mm

\section{Introduction}
To study the local behavior of solutions to second order elliptical partial differential equations, Morrey [12] introduced  Morrey space,  which is defined as
$$M^{p,\lambda}(\Bbb{R}^n)=\left\{f\in M^{p,\lambda}(\Bbb{R}^n):\|f\|_{M^{p,\lambda}(\Bbb{R}^n)}:=\sup\limits_B\frac{1}{|B|^\lambda}\left(\frac{1}{|B|}\int_B|f (x)|^pdx\right)^{1/p}<\infty\right\}$$
with the exponents $p$ and $\lambda$ satisfying  $p\geq 1$ and $-\frac{1}{p}<\lambda<0$.

  For the extension of Morrey space $M^{p,\lambda}(\Bbb{R}^n)$, the classical Campanato space $\mathcal{C}^{p,\lambda}(\Bbb{R}^n)$ is defined by

  $$\mathcal{C}^{p,\lambda}(\Bbb{R}^n)=\left\{f\in \mathcal{C}^{p,\lambda}(\Bbb{R}^n):\|f\|_{\mathcal{C}^{p,\lambda}(\Bbb{R}^n)}:=\sup\limits_B\frac{1}{|B|^\lambda}\left(\frac{1}{|B|}\int_B|f(x)-f_B|^pdx\right)^{1/p}<\infty\right\},$$
  where $p\in [1,\infty)$, $\lambda\in (-\frac{1}{p}, \frac{1}{n})$, $f_B=\frac{1}{|B|}\int_Bf(x)dx$ and $B\subset \Bbb{R}^n$ denotes any ball in $\Bbb{R}^n$.

  It is well known that $M^{p,\lambda}(\Bbb{R}^n)\subset \mathcal{C}^{p,\lambda}(\Bbb{R}^n)$. For the studies of such two function spaces and the action of various operators on them, one may see [14-16] et al. for more details.

  In  [10, 11], Lu and Yang studied  a new kind of  homogeneous  Hardy type space $H\dot{A}_q$ with $q>1$ and  they found that the dual space of $H\dot{A}_q$ can be defined by the following norm.
  $$
  \parallel f\parallel_{C\dot{B}MO^q}:=\sup\limits_{R\geq
  0}\left(\frac{1}{|B(0,R)|}\int_{B(0,R)}|f(x)-f_{B(0,R)}|^qdx\right)^{1/q}<\infty.
  $$

  Obviously,  ${C\dot{B}MO}^q$ is the homogeneous central bounded mean oscillation depending on $q$. Moreover, the famous John-Nirenberg inequality no longer holds in such space. Thus, it can be regarded as an extension of the classical BMO space.

   In 2000,  Alvarez, Lakey and Guzm\'{a}n-Partida [1] introduced the  $\lambda$-central Campanato space and $\lambda$-central Morrey space respectively.

\begin{definition} ([1]) Let $-\frac{1}{p}<\lambda<\frac{1}{n}$ with $1<p<\infty$. Then, a function $f\in L^p_{\text{loc}}(\Bbb{R}^n)$ is said to  belonged to the $\lambda$-central Campanato space
  $\dot{\mathcal{C}}^{p,\lambda}(\Bbb{R}^n)$ if
   $$
  \parallel
  f\parallel_{\dot{\mathcal{C}}^{p,\lambda}(\Bbb{R}^n)}:=\sup\limits_{R>0}\left(\frac{1}{|B(0,R)|^{1+\lambda
  p}}\int_{B(0,R)}|f(x)-f_{B(0,R)}|^pdx\right)^{1/p}<\infty.$$
\end{definition}
\begin{definition} ([1])  Let $-\frac{1}{p}<\lambda<\frac{1}{n}$ and
  $1<p<\infty$. Then, the $\lambda$-central Morrey space
  $\dot{\mathcal{M}}^{p,\lambda}(\Bbb{R}^n)$ is defined by

   $$
  \dot{\mathcal{M}}^{p,\lambda}(\Bbb{R}^n)=\left\{  \parallel
  f\parallel_{\dot{\mathcal{M}}^{p,\lambda}(\Bbb{R}^n)}=\sup\limits_{R>0}\left(\frac{1}{|B(0,R)|^{1+\lambda
  p}}\int_{B(0,R)}|f(x)|^pdx\right)^{1/p}<\infty\right\}.
 $$\end{definition}

   If $0<\lambda<1/p$,  the $\lambda$-central-Campanato space becomes the $\lambda$-central bounded mean oscillation space $C\dot{B}MO^{p,\lambda}(\Bbb{R}^n)$. Moreover, it is easy to check that $\mathcal{C}^{p,\lambda}(\Bbb{R}^n)\subset \dot{\mathcal{C}}^{p,\lambda}(\Bbb{R}^n)$ and $\dot{\mathcal{C}}^{p,0}(\Bbb{R}^n)=C\dot{{B}}MO^p(\Bbb{R}^n)$.   For the case $-1/p<\lambda<0$, there is $\dot{\mathcal{M}}^{p,\lambda}(\Bbb{R}^n)\subset \dot{\mathcal{C}}^{p,\lambda}(\Bbb{R}^n)$.

Suppose that  $T$ is an integral operator and $b$ is a local integrable function. Then,  the commutator of $T$ is defined by

   $$T_b(f)(x)=b(x)T(f)(x)-T(b)(x).$$

   For the actions  of commutators on $\dot{\mathcal{M}}^{p,\lambda}(\Bbb{R}^n)$,  one may see [4,8,17] et al. to find more details  with $b\in \dot{\mathcal{C}}^{p,\lambda}(\Bbb{R}^n)$ and $0<\lambda<\frac{1}{n}$.  For the case $b\in \dot{\mathcal{C}}^{p,\lambda}(\Bbb{R}^n)$ with $-\frac{1}{p}<\lambda<0$, Shi and Lu [16] studied the boundedness of commutators generated by the Hardy operators.

   Next, we give some definitions about Hardy operators.

   For  $f\in L^p(\Bbb{R}^+)$ with $1<p<\infty$,  the classical Hardy operator is  defined by
   $$
   Hf(x)=\frac{1}{x}\int_0^xf(t)dt, x\neq 0.
      $$

   In 1920,  Hardy [7] proved the $L^p(\Bbb{R}^+)$ boundedness of $H$ and showed the constant $\frac{p}{p-1}$ of  (1) is the best possible.
   \begin{equation}\label{resolution}
    \parallel Hf\parallel_{L^p(\Bbb{R}^+)}\leq \frac{p}{p-1}\parallel f\parallel_{L^p(\Bbb{R}^+)}.
   \end{equation}

In 1995, Christ and Grafakos [2] introduced  the following $n$-dimensional Hardy operator $\mathcal{H}$ defined by
$$
\mathcal{H}f(x)=\frac{1}{|x|^{n}}\int_{|y|<|x|}f(y)dy, \ \ \ \ \ \ \ \ \ x\in \Bbb{R}^n\setminus \{0\},
  $$
  Christ and Grafakos [2] also showed  the operator $\mathcal{H}$ satisfies the analogue results of (1).

The dual operator of $\mathcal{H}$ is $\mathcal{H}^\ast$, which is defined by
$$
\mathcal{H}^\ast f(x)=\int_{|y|\geq |x|}\frac{f(y)}{|y|^n}dy, \ \ \ \ \ \ \ \ \ x\in \Bbb{R}^n\setminus \{0\}.
  $$

It is easy to check that
$$\int_{\Bbb{R}^n}g(x)\mathcal{H}f(x)dx=\int_{\Bbb{R}^n}f(x)\mathcal{H}^\ast g(x)dx.$$

Thus, the commutator of  Hardy type operator is defined as
$$
\mathcal{H}_bf(x)=b(x)\mathcal{H}f(x)-\mathcal{H}(fb)(x) \ \ \ \text{and} \ \ \ \mathcal{H}^\ast_bf(x)=b(x)\mathcal{H}^\ast f(x)-\mathcal{H}^\ast (fb)(x).
$$

The operators $\mathcal{H}_b$ and $\mathcal{H}^\ast_b$  were first studied in [5]  where Fu et al. gave the characterization of $C\dot{B}MO^p(\Bbb{R}^n)$ with $1<p<\infty$ via the boundedness of $\mathcal{H}_b$ and $\mathcal{H}^\ast_b$  on $L^p(\Bbb{R}^n)$. For more studies about the operartor $\mathcal{H}_b$ and $\mathcal{H}^\ast_b$, one may see the paper [16,19] to find more details about $\mathcal{H}_b$ and $\mathcal{H}^\ast_b$ with $b\in  \dot{\mathcal{C}}^{p,\lambda}(\Bbb{R}^n)$.

On the other hand, the weighted norm inequalities for integral operators was first studied in the last 70s  and  one may  see [3,13] et. al. for more details. In 2009, Komori and Shirai [9]  defined the weighted Morrey space and they showed the boundedness of some classical integral operators and their commutators on the weighted Morrey spaces.

    Next, we  introduce the weighted central-Campanato space $\dot{\mathcal{C}}^{p,\lambda}_\omega(\Bbb{R}^n)$ and weighted $\lambda$-central Morrey space
$\dot{\mathcal{M}}^{p,\lambda}_{\omega}(\Bbb{R}^n)$ respectively by the following norms.

$$\|f\|_{\dot{\mathcal{M}}^{p,\lambda}_\omega(\Bbb{R}^n)} :=\sup\limits_{R>0}\left(\frac{1}{\omega(B(0,R)^{1+\lambda
q}}\int_{B(0,R)}|f(x)|^p\omega(x)dx\right)^{1/p}$$
and

$$\|f\|_{\dot{\mathcal{C}}^{p,\lambda}_\omega(\Bbb{R}^n)}:=\sup\limits_{R>0}\left(\frac{1}{\omega(B(0,R))^{1+\lambda p}}\int_{B(0,R)}|f(x)- f_{B,\omega(B(0,R))}|^p\omega(x)dx\right)^{1/p}$$
where the definition of $f_{B,\omega}$ is  $f_{B,\omega}=\frac{1}{\omega(B)}\int_{B}f(x)\omega(x)dx$ and the exponents of $p,\lambda$ are the same as in the definition of  ${\mathcal{C}}^{p,\lambda}(\Bbb{R}^n)$ and ${\mathcal{M}}^{p,\lambda}(\Bbb{R}^n)$.

For the boundedness of integral operators on $\dot{\mathcal{M}}^{p,\lambda}_\omega(\Bbb{R}^n)$ with $\lambda<0$, one may see [18] et al. for more details.

 Suppose that  $\omega$ is a non-negative and locally integrable function. If for every cube $Q$, there exists a constant $D$ independent of $Q$, such that $\omega(2Q)\leq D\omega(Q)$. Then, we say $\omega$ satisfy the doubling  condition and we simply denote $\omega\in \Delta_2$.

We would like to mention that in [18], the restriction of $\omega$ is $\omega\in A_p$ where $A_p$ denotes the Muckenhoupt weight classes (see [13]). From [6], we know that if $\omega\in A_p$, then $\omega$ satisfies the doubling condition (i.e. $\omega\in \Delta_2$). However, the converse is not true.  Throughout this paper,  we only assume that $\omega\in \Delta_2.$

Now, we  are interested in the following weighted Hardy operator  $\mathcal{H}_\omega$ with the weight  $\omega\in \Delta_2$.
$$\mathcal{H}_\omega(f)(x)=\frac{1}{\omega(B(0,|x|))}\int_{|y|<|x|}f(y)\omega(y)dy.$$

For any $g\in L^1_{\text{loc}}(\omega)$ with $\omega\in \Delta_2$, we have
\begin{align*}
<\mathcal{H}_\omega(f)(x),g>_\omega&=\int_{\Bbb{R}^n}\frac{1}{\omega(B(0,|x|))}\int_{|y|<|x|}f(y)\omega(y)dyg(x)\omega(x)dx\\
&=\int_{\Bbb{R}^n}\int_{|x|>|y|}\frac{g(x)\omega(x)}{\omega(B(0,|x|))}dxf(y)\omega(y)dy.
\end{align*}

Thus, we can define the dual operator of $\mathcal{H}_\omega$ as
$$\mathcal{H}^\ast_\omega(g)(x)=\int_{|y|\geq |x|}\frac{g(y)\omega(y)}{\omega(B(0,|y|))}dy.$$

Then, the commutators of $\mathcal{H}_{\omega,b}$ and $\mathcal{H}_{\omega,b}^\ast$ can be stated as follows.
$$\mathcal{H}_{\omega,b}(f)(x)=b(x)\mathcal{H}_{\omega}(f)(x)-\mathcal{H}_{\omega}(fb)(x)$$
and
$$\mathcal{H}^\ast_{\omega,b}(f)(x)=b(x)\mathcal{H}^\ast_{\omega}(f)(x)-\mathcal{H}^\ast_{\omega}(fb)(x).$$

In this paper, we will give the boundedness  of $\mathcal{H}_{\omega,b}$ and $\mathcal{H}_{\omega,b}^\ast$ on $\dot{\mathcal{M}}^{p,\lambda}_\omega(\Bbb{R}^n)$ with $\dot{\mathcal{C}}^{p,\lambda}_\omega(\Bbb{R}^n)$ and  $\lambda<0$.

\begin{theorem} Let $1<p<\infty$, $-\frac{1}{p}<\lambda<0$, $-\frac{1}{p_i}<\lambda_i<0$ with $i=1,2$, $\frac{1}{p}=\sum\limits_{i=1}^2\frac{1}{p_i}$ and  $\lambda=\sum\limits_{i=1}^2\lambda_i$. Then,   both $\mathcal{H}_{\omega,b}$ and $\mathcal{H}^\ast_{\omega,b}$ are bounded from $\dot{\mathcal{M}}^{p_2,\lambda_2}_\omega(\Bbb{R}^n)$ to $\dot{\mathcal{M}}^{p,\lambda}_\omega(\Bbb{R}^n)$ with $b\in \dot{\mathcal{C}}^{p_1,\lambda_1}_\omega(\Bbb{R}^n)$.
\end{theorem}

Morover, we have

\begin{theorem} Let $2<p<\infty$ and  $-\frac{1}{2p}<\lambda<0$. Then,   Both $\mathcal{H}_{\omega,b}$ and $\mathcal{H}^\ast_{\omega,b}$ are bounded from $\dot{\mathcal{M}}^{p,\lambda}_\omega(\Bbb{R}^n)$ to $\dot{\mathcal{M}}^{p,2\lambda}_\omega(\Bbb{R}^n)$ with $ b\in \dot{\mathcal{C}}^{p,\lambda}_\omega(\Bbb{R}^n)$.
\end{theorem}
In order to give the characterization of $\dot{\mathcal{C}}^{p,\lambda}_\omega(\Bbb{R}^n)$, we introduce the  operators $\mathcal{H}_{\omega,|b|}$ and $\mathcal{H}^\ast_{\omega,|b|}$  as follows.

$$\mathcal{H}_{\omega,|b|}=\frac{1}{\omega(B(0,|x|))}\int_{|y|<|x|}f(y)|b(x)-b(y)|\omega(y)dy$$
and
$$\mathcal{H}^\ast_{\omega,|b|}(f)(x)=\int_{|y|\geq |x|}\frac{f(y)\omega(y)}{\omega(B(0,|y|))}|b(x)-b(y)|dy.$$

By checking the proofs of Theorems 1.1 and 1.2, we know that the above two theorems still hold if we replace $\mathcal{H}_{\omega,b}$ by $\mathcal{H}_{\omega,|b|}$ and $\mathcal{H}^\ast_{\omega,b}$ by $\mathcal{H}^\ast_{\omega,|b|}$. Moreover, we have the following theorems.

\begin{theorem} Let $1<p<\infty$, $-\frac{1}{p}<\lambda<0$, $-\frac{1}{p_i}<\lambda_i<0$ with $i=1,2$, $\frac{1}{p}=\sum\limits_{i=1}^2\frac{1}{p_i}$ and  $\lambda=\sum\limits_{i=1}^2\lambda_i$.  Moreover, we assume that $b$ satisfies
\begin{equation}\label{resolution}
\sup\limits_{B(0,R)\ni x}|b(x)-b_{B,\omega{\left(B(0,R)\right)}}|\leq \frac{C}{\omega(B(0,R))}\int_{B(0,R)}|b(x)-b_{B,\omega{\left(B(0,R)\right)}}|\omega(x)dx, \end{equation}
Then,  the following two conditions are equivalent.

 ($a_3$) Both $\mathcal{H}_{\omega,|b|}$ and $\mathcal{H}^\ast_{\omega,|b|}$ are bounded from $\dot{\mathcal{M}}^{p_2,\lambda_2}_\omega(\Bbb{R}^n)$ to $\dot{\mathcal{M}}^{p,\lambda}_\omega(\Bbb{R}^n)$.

($b_3$) $b\in \dot{\mathcal{C}}^{p_1,\lambda_1}_\omega(\Bbb{R}^n)$.\\
\end{theorem}
In order to give up the condition (2), we have

\begin{theorem}
  Let $2<p<\infty$ and  $-\frac{1}{2p}<\lambda<0$. Then,   the following two conditions are equivalent.

 ($a_4$) Both $\mathcal{H}_{\omega,b}$ and $\mathcal{H}^\ast_{\omega,b}$ are bounded from $\dot{\mathcal{M}}^{p,\lambda}_\omega(\Bbb{R}^n)$ to $\dot{\mathcal{M}}^{p,2\lambda}_\omega(\Bbb{R}^n)$.

 ($b_4$) $ b\in \dot{\mathcal{C}}^{p,\lambda}_\omega(\Bbb{R}^n)$.

\end{theorem}

\section{Some useful lemmas.}
 In this section, we give some lemmas  that will be used throughout this paper.  For simplicity, we denote $B=B(0,R)$,  $B_i=\{x\in \Bbb{R}^n:|x|\leq 2^i\}, \ \ C_i=B_i\setminus B_{i-1}$ with  $i\in \Bbb{Z}$  and $C$ may represents different constants in different places.

\begin{lemma} ([6])  If $\omega\in \Delta_2$, then there exists a constant $D_1: D_1>1$ independent of $Q$ such that
$\omega(2Q)\leq {D_1}\omega(Q),$ for any cube $Q$.\\
\end{lemma}
\begin{lemma} ([9, Lemma 4.1])  If $\omega\in \Delta_2$, then there exists a constant $D_2:D_2>1$  independent of $Q$ such that
$\omega(2Q)\geq {D_2}\omega(Q)$. \\
\end{lemma}
\begin{rem}  By checking the proof of [9, Lemma 4.1], we know that the constant $D_2$ is strictly less than 2. Moreover, there exists a constant $r>1$ independent of $Q$, such that $r\leq D_2<2$.
\end{rem}

\begin{lemma} \  For  $\forall i,\ j \in \Bbb{Z}$ and $\forall\omega\in \Delta_2$, we have the following inequalities.

(i) If $i\geq j$, there is ${D_2}^{i-j}\omega\left(B_{j}\right)\leq \omega\left(B_i\right)\leq{D_1}^{i-j}\omega\left(B_{j}\right)$.

(ii) For the case $i\leq j$, we have ${D_1}^{i-j}\omega\left(B_{j}\right)\leq \omega\left(B_i\right)\leq{D_2}^{i-j}\omega\left(B_{j}\right)$.

Proof.  Lemma 2.3 is a simple derivation of Lemmas  2.1  and    2.2 and we omit the proof process.
\end{lemma}

\begin{lemma}  For $\forall\omega\in \Delta_2$ and  $x\in C_i$ with $i\in \Bbb{Z}^+$, there is $\frac{1}{D_1}\omega(B_i)\leq\omega(B(0,|x|)\leq\omega(B_i)$.

 Proof. As  $i\in \Bbb{Z}$ and $x\in C_i$, it is easy to see $\omega(B_{i-1})\leq\omega(B(0,|x|)\leq\omega(B_i)$.
    Since $\omega\in \Delta_2$, we get $\frac{1}{D_1}\omega(B_i)\leq\omega(B_{i-1})$.

    Thus, we obtain
    $$\frac{1}{D_1}\omega(B_i)\leq\omega(B(0,|x|)\leq\omega(B_i).$$

\end{lemma}

\begin{lemma} \  Let $1<p<\infty$, $-\frac{1}{p}<\lambda<0$ and  $b\in \dot{\mathcal{C}}^{p,\lambda}_\omega(\Bbb{R}^n)$. Then,  for any $i,\ j \in \Bbb{Z}$ with $i<j$ and $\forall\omega\in \Delta_2$, we have
$$|b_{B_{i,\omega}}-b_{B_j,\omega}|\leq C \|b\|_{\dot{\mathcal{C}}^{p,\lambda}_\omega(\Bbb{R}^n)}\omega(B_i)^\lambda.$$ \\

Proof. As $i<j$, using Lemma 2.3 and the H\"{o}lder inequality with $1<p<\infty$ and  $-\frac{1}{p}<\lambda<0$, we obtain
\begin{align*}
|b_{B_{i},\omega}-&b_{B_{j},\omega}|\leq\sum_{k=i}^{j-1}\left|b_{B_k,\omega}-b_{B_{k+1},\omega}\right|\\
&\leq\sum_{k=i}^{j-1}\frac{1}{\omega(B_{k})}\int_{B_{k}}\left|b(x)-b_{B_{k+1},\omega}\right|\omega(x)dx\\
&\leq\sum_{k=i}^{j-1}\frac{1}{\omega(B_{k})}\left(\int_{B_{k+1}}\left|b(x)-b_{B_{k+1},\omega}\right|^{p}\omega(x)dx\right)^{\frac{1}{p}}\left(\int_{B_{k+1}}\omega(x)dx\right)^{\frac{1}{p'}}\\
&\leq C\sum_{k=i}^{j-1}\frac{1}{\omega(B_{k})}\left(\omega(B_{k+1})\right)^{\lambda+\frac{1}{p}}{||b||}_{\dot{\mathcal{C}}^{p,\lambda}_\omega(\Bbb{R}^n)
  }\left(\omega(B_{k+1})\right)^{\frac{1}{p'}}\\
&\leq CD_1{||b||}_{\dot{\mathcal{C}}^{p,\lambda}_\omega(\Bbb{R}^n)}{\omega(B_{i})}^{\lambda}\sum_{k=i}^{j-1}{D_2}^{(k+1-i)\lambda}\\
&\leq CD_1\|b\|_{\dot{\mathcal{C}}^{p,\lambda}_\omega(\Bbb{R}^n)}\omega(B_{i})^{\lambda},\\
\end{align*}
here $C$  is a positive constant independent of $j$ and $i$.
\end{lemma}

\begin{lemma}\  Let  $1<p<\infty$, $-\frac{1}{p}<\lambda<0$ and  $f\in L^p_\omega(\Bbb{R}^n)$ with $\omega\in \Delta_2$.  If $\int_{B_{k}}|f(x)|^p\omega(x)dx\leq C\omega(B_{k})^{1+\lambda p}$ holds for any $k\in \Bbb{Z}$,
 then  for  any $R>0$,  we have
 $$\int_{B{(0,R)}}|f(x)|^p\omega(x)dx\leq C\omega(B{(0,R)})^{\lambda p+1}.$$

Proof. For any $R>0$, there exists a $k\in \Bbb{Z}$, such that $2^{k-1}<R\leq 2^{k}$.  Thus,  we get $B_{k-1}\subset B{(0,R)}\subset B_{k}$ and $\omega\left(B_{k-1}\right)\leq\omega\left(B{(0,R)}\right)\leq\omega\left(B_{k}\right)$.  Then,  using  Lemmas 2.3- 2.4, we get

 \begin{align*}
&\int_{B{(0,R)}}|f(x)|^p\omega(x)dx
\leq\int_{B_{k}}|f(x)|^p\omega(x)dx\\
&\leq C\omega(B_{k})^{1+\lambda p}
 \leq C\left(D_1\omega(B_{k-1})\right)^{\lambda p+1}\\
&\leq C\omega(B{(0,R)})^{\lambda p+1}.
\end{align*}
\end{lemma}

\begin{lemma}  Let $1<p<\infty$ and $-\frac{1}{p}<\lambda<0$. Then for  $\forall R> 0$ and  $\forall\omega\in \Delta_2$, there is $$\|\chi_B\|_{\dot{\mathcal{M}}^{p,\lambda}_\omega(\Bbb{R}^n)}\leq {\omega(B)^{-\lambda}}.$$

Proof. By the definition of $\dot{\mathcal{M}}^{p,\lambda}_\omega(\Bbb{R}^n)$, we have
\begin{align*}
\|\chi_B\|_{\dot{\mathcal{M}}^{p,\lambda}_\omega(\Bbb{R}^n)}
&=\sup\limits_{r>0}\left(\frac{1}{\left(\omega(B(0,r)\right)^{1+\lambda p}}\int_{B(0,r)}|\chi_B(x)|^p\omega(x)dx\right)^{1/p}\\
&=\sup\limits_{r>0}\left(\frac{1}{\left(\omega(B(0,r)\right)^{1+\lambda p}}\int_{B(0,r)\cap B(0,R)}\omega(x)dx\right)^{1/p},\\
\end{align*}

For the case $r\leq R$, as since $1<p<\infty$ and $-\frac{1}{p}<\lambda<0$, there is
\begin{align*}
&\left(\frac{1}{\left(\omega(B(0,r)\right)^{1+\lambda p}}\int_{B(0,r)\cap B(0,R)}\omega(x)dx\right)^{1/p}=\left(\frac{1}{\left(\omega(B(0,r)\right)^{1+\lambda p}}\int_{B(0,r)}\omega(x)dx\right)^{1/p}\\
&=\left(\frac{1}{\omega(B(0,r)^{\lambda p}}\right)^{1/p}\leq\left(\frac{1}{\omega(B(0,R)^{\lambda  p}}\right)^{1/p}={\omega(B)^{-\lambda}}, \\
\end{align*}

For the case $r>R$, we get
\begin{align*}
&\left(\frac{1}{\left(\omega(B(0,r)\right)^{1+\lambda p}}\int_{B(0,r)\cap B(0,R)}\omega(x)dx\right)^{1/p}=\left(\frac{1}{\left(\omega(B(0,r)\right)^{1+\lambda p}}\int_{B(0,R)}\omega(x)dx\right)^{1/p}\\
&\leq\left(\frac{1}{\left(\omega(B(0,R)\right)^{1+\lambda p}}\int_{B(0,R)}\omega(x)dx\right)^{1/p}=\left(\frac{1}{\left(\omega(B(0,R)\right)^{\lambda p}}\right)^{1/p}= {\omega(B)^{-\lambda}},\\
\end{align*}
Thus, for any $r>0$, we have

$$\sup\limits_{r>0}\left(\frac{1}{\left(\omega(B(0,r)\right)^{1+\lambda p}}\int_{B(0,r)\cap B(0,R)}\omega(x)dx\right)^{1/p}\leq {\omega^{- \lambda}(B)}, $$
which finishes the proof of Lemma 2.7.
\end{lemma}

\section{Proof of Theorem 1.1.}

By the definitions of  $\dot{\mathcal{M}}^{p,\lambda}_\omega(\Bbb{R}^n)$ and Lemma 2.6,  it suffices to show   the following estimates with  $k\in \Bbb{Z}$.

\begin{equation}\label{resolution}\int_{B_{k}}|\mathcal{H}_{\omega,b}(f)(x)|^p\omega(x)dx\leq C\omega(B_{k})^{1+\lambda p}\|b\|_{\dot{\mathcal{C}}^{p_1,\lambda_1}_\omega(\Bbb{R}^n)}^p\|f\|^p_{\dot{\mathcal{M}}^{p_2,\lambda_2}_\omega(\Bbb{R}^n)}\end{equation}

\begin{equation}\label{resolution}\int_{B_{k}}|\mathcal{H}^\ast_{\omega,b}(f)(x)|^p\omega(x)dx\leq C\omega(B_{k})^{1+\lambda p}\|b\|_{\dot{\mathcal{C}}^{p_1,\lambda_1}_\omega(\Bbb{R}^n)}^p\|f\|^p_{\dot{\mathcal{M}}^{p_2,\lambda_2}_\omega(\Bbb{R}^n)}.\end{equation}

We begin with the proof of (3).  Using Lemma 2.4, we get

\begin{align*}
&\int_{B_{k}}\left(\mathcal{H}_{\omega,b}(f)(x)\right)^p\omega(x)dx\\
&\leq \sum\limits_{j=-\infty}^{k}\int_{C_j}\left(\frac{1}{\omega(B(0,|x|)}\int_{B(0,|x|)}|b(x)-b(y)||f(y)|\omega(y)dy\right)^p\omega(x)dx\\
&\leq \sum\limits_{j=-\infty}^{k}\frac{1}{\omega(B_j)^p}\int_{B_j}\left(\int_{B_j}|b(x)-b(y)||f(y)|\omega(y)dy\right)^p\omega(x)dx\\
&\leq \sum\limits_{j=-\infty}^{k}\frac{1}{\omega(B_j)^p}\int_{B_j}\left(\sum\limits_{i=-\infty}^j\int_{B_i}|b(x)-b(y)||f(y)|\omega(y)dy\right)^p\omega(x)dx\\
&\leq \sum\limits_{j=-\infty}^{k}\frac{1}{\omega(B_j)^p}\int_{B_j}\left(\sum\limits_{i=-\infty}^j\int_{B_i}|b(x)-b_{B_j,\omega}||f(y)|\omega(y)|dy\right)^p\omega(x)dx\\
& \ \ \ +\sum\limits_{j=-\infty}^{k}\frac{1}{\omega(B_j)^p}\int_{B_j}\left(\sum\limits_{i=-\infty}^j\int_{B_i}|b(y)-b_{B_j,\omega}||f(y)|\omega(y)|dy\right)^p\omega(x)dx\\
&=:I_1+I_2.
\end{align*}

For $I_1$,  we prove the fact.
\begin{equation}\label{resolution}
\sum\limits_{i=-\infty}^j\int_{B_i}|f(y)|\omega(y)dy
\leq C\|f\|_{\dot{\mathcal{M}}^{p_2,\lambda_2}_\omega(\Bbb{R}^n)}\left(\omega(B_j)\right)^{\lambda_2+1} \end{equation}
 with  $1<p_2<\infty$ and $-1<\lambda_2<0$.

By Lemma 2.3 and the H\"{o}lder inequality, there is

\begin{align*}
&\sum\limits_{i=-\infty}^j\int_{B_i}|f(y)|\omega(y)dy\\
& \leq C\sum\limits_{i=-\infty}^j\left(\int_{B_i}\left(|f(y)|\omega^\frac{1}{p_2}(y)\right)^ {p_2}  dy\right)^\frac{1}{p_2}\left(\int_{B_i}\omega(y)^{(1-\frac{1}{p_2})p_2'}dy\right)^{1/\left(\frac{1}{p_2}\right)'}\\
&\leq C\|f\|_{\dot{\mathcal{M}}^{p_2,\lambda_2}_\omega(\Bbb{R}^n)}\sum\limits_{i=-\infty}^j\left(\omega(B_i)\right)^{\lambda_2+\frac{1}{p_2}}\left(\omega(B_i)\right)^{\left(1-\frac{1}{p_2}\right)}\\
&\leq C\|f\|_{\dot{\mathcal{M}}^{p_2,\lambda_2}_\omega(\Bbb{R}^n)}\sum\limits_{i=-\infty}^j\left(\omega(B_i)\right)^{\lambda_2+1}\\
&\leq C\|f\|_{\dot{\mathcal{M}}^{p_2,\lambda_2}_\omega(\Bbb{R}^n)}\sum\limits_{i=-\infty}^j\left({D_2}^{(i-j)}\omega(B_j)\right)^{\lambda_2+1}\\
&\leq C\|f\|_{\dot{\mathcal{M}}^{p_2,\lambda_2}_\omega(\Bbb{R}^n)}\left(\omega(B_j)\right)^{\lambda_2+1}\\
\end{align*}
and we finish the proof of (5).

Then,  using  Lemma 2.3, the H\"{o}lder inequality and the conditions of Theorem 1.1,  we get

\begin{align*}
I_1
&=C\sum\limits_{j=-\infty}^{k}\frac{1}{\omega(B_j)^p}\int_{B_j}|b(x)-b_{B_j,\omega}|^p\omega(x)dx\left(\sum\limits_{i=-\infty}^j\int_{B_i}|f(y)|\omega(y)dy\right)^p\\
&\leq C\sum\limits_{j=-\infty}^{k}\frac{1}{\omega(B_j)^p}\left(\int_{B_j}\left(|b(x)-b_{B_j,\omega}|^p\omega^\frac{p}{p_1}(x)\right)^\frac{p_1}{p}dx\right)^\frac{p}{p_1}\left(\int_{B_j}\omega(x)^{(1-\frac{p}{p_1})\left(\frac{p_1}{p}\right)'}dx\right)^{1/\left(\frac{p_1}{p}\right)'}\\
& \ \ \ \times \left(\|f\|_{\dot{\mathcal{M}}^{p_2,\lambda_2}_\omega(\Bbb{R}^n)}\left(\omega(B_j)\right)^{\lambda_2+1}\right)^p\\
&\leq C\sum\limits_{j=-\infty}^{k}\frac{1}{\omega(B_j)^p}\left(\left(\omega(B_j)\right)^{{\lambda_1+\frac{1}{p_1}}}\|b\|_{\dot{\mathcal{C}}^{p_1,\lambda_1}_\omega(\Bbb{R}^n)}\right)^p\left(\omega(B_j)\right)^{1/\left(\frac{p_1}{p}\right)'}\\
&\times\left(\|f\|_{\dot{\mathcal{M}}^{p_2,\lambda_2}_\omega(\Bbb{R}^n)}\left(\omega(B_j)\right)^{\lambda_2+1}\right)^p\\
&\leq C\|b\|_{\dot{\mathcal{C}}^{p_1,\lambda_1}_\omega(\Bbb{R}^n)}^p\|f\|_{\dot{\mathcal{M}}^{p_2,\lambda_2}_\omega(\Bbb{R}^n)}^p
  \sum\limits_{j=-\infty}^{k}\left(\omega(B_j)\right)^{{1+\lambda p}} \\
&\leq C\|b\|_{\dot{\mathcal{C}}^{p_1,\lambda_1}_\omega(\Bbb{R}^n)}^p\|f\|_{\dot{\mathcal{M}}^{p_2,\lambda_2}_\omega(\Bbb{R}^n)}^p
  \sum\limits_{j=-\infty}^{k}\left({D_2}^{\left(j-k\right)}\omega(B_{k})\right)^{{1+\lambda p}}\\
&\leq C\|b\|_{\dot{\mathcal{C}}^{p_1,\lambda_1}_\omega(\Bbb{R}^n)}^p\|f\|_{\dot{\mathcal{M}}^{p_2,\lambda_2}_\omega(\Bbb{R}^n)}^p
  \left(\omega(B_{k})\right)^{{1+\lambda p}}.\\
\end{align*}

For  $I_2$,   we need to show the following estimates.
\begin{equation}\label{resolution}
\sum\limits_{i=-\infty}^j\int_{B_i}|b(y)-b_{B_j,\omega}||f(y)|\omega(y)dy
\leq C\|b\|_{\dot{\mathcal{C}}^{p_1,\lambda_1}_\omega(\Bbb{R}^n)}\|f\|_{\dot{\mathcal{M}}^{p_2,\lambda_2}_\omega(\Bbb{R}^n)}\left(\omega(B_j)\right)^{1+\lambda}. \end{equation}
where $x\in C_j$, $1<p_i<\infty$, $-\frac{1}{p_i}<\lambda_i<0$ with $i=1,2$ and  $\lambda=\sum\limits_{i=1}^2\lambda_i$.

Denote $s$ by $1/s=1-1/{p_1}-1/{p_2}$.  Using Lemma 2.3 and the H\"{o}lder inequality, we get

\begin{align*}
&\sum\limits_{i=-\infty}^j\int_{B_i}|b(y)-b_{B_j,\omega}||f(y)|\omega(y)dy\\
&\leq C\sum\limits_{i=-\infty}^j\left(\int_{B_i}|b(y)-b_{B_j,\omega}|^{p_1}\omega(y)dy\right)^{1/p_1}\left(\int_{B_i}|f(y)|^{p_2}\omega(y)dy\right)^{1/p_2}\left(\int_{B_i}\omega(y)^{\frac{1}{s}s}dy\right)^{1/s}\\
&\leq C\|b\|_{\dot{\mathcal{C}}^{p_1,\lambda_1}_\omega(\Bbb{R}^n)}\|f\|_{\dot{\mathcal{M}}^{p_2,\lambda_2}_\omega(\Bbb{R}^n)}\sum\limits_{i=-\infty}^j\omega(B_i)^{1+\lambda}\\
&\leq C\|b\|_{\dot{\mathcal{C}}^{p_1,\lambda_1}_\omega(\Bbb{R}^n)}\|f\|_{\dot{\mathcal{M}}^{p_2,\lambda_2}_\omega(\Bbb{R}^n)} \sum\limits_{i=-\infty}^j\left({D_2}^{(i-j)}\omega(B_j)\right)^{1+\lambda}\\
&\leq C\|b\|_{\dot{\mathcal{C}}^{p_1,\lambda_1}_\omega(\Bbb{R}^n)}\|f\|_{\dot{\mathcal{M}}^{p_2,\lambda_2}_\omega(\Bbb{R}^n)}\omega(B_{j})^{1+\lambda}\\
\end{align*}
and we finish the proof of (6).

Thus, we obtain

\begin{align*}
I_2&\leq C\sum\limits_{j=-\infty}^{k}\int_{B_j}\frac{1}{\omega(B_j)^p}\left(\sum\limits_{i=-\infty}^j\int_{B_i}|b(y)-b_{B_j,\omega}||f(y)|\omega(y)dy\right)^p\omega(x)dx\\
&\leq C\|b\|^p_{\dot{\mathcal{C}}^{p_1,\lambda_1}_\omega(\Bbb{R}^n)}\|f\|^p_{\dot{\mathcal{M}}^{p_2,\lambda_2}_\omega(\Bbb{R}^n)}\sum\limits_{j=-\infty}^{k}\int_{B_j}\frac{1}{\omega(B_j)^p}\left(\omega(B_j)^{1+\lambda}\right)^p\omega(x)dx\\
&\leq \|b\|^p_{\dot{\mathcal{C}}^{p_1,\lambda_1}_\omega(\Bbb{R}^n)}\|f\|^p_{\dot{\mathcal{M}}^{p_2,\lambda_2}_\omega(\Bbb{R}^n)}\sum\limits_{j=-\infty}^{k}\omega(B_j)^{1+\lambda p}\\
&\leq \|b\|^p_{\dot{\mathcal{C}}^{p_1,\lambda_1}_\omega(\Bbb{R}^n)}\|f\|^p_{\dot{\mathcal{M}}^{p_2,\lambda_2}_\omega(\Bbb{R}^n)}\omega(B_{k})^{1+\lambda p}\sum\limits_{j=-\infty}^{k}{D_2}^{(j-k)(1+\lambda p)}\\
&\leq C\|b\|^p_{\dot{\mathcal{C}}^{p_1,\lambda_1}_\omega(\Bbb{R}^n)}\|f\|^p_{\dot{\mathcal{M}}^{p_2,\lambda_2}_\omega(\Bbb{R}^n)}\omega(B_{k})^{1+\lambda p}.
\end{align*}
Combing the estimates of $I_1$ and $I_2$, we finish the proof of  (3).

Next, we will prove (4). First,  we can decompose $\int_{B_{k}}|\mathcal{H}^\ast_{\omega,b}(f)(x)|^p\omega(x)dx$ as follows.

 \begin{align*}
\int_{B_{k}}|\mathcal{H}^\ast_{\omega,b}(f)(x)|^p\omega(x)dx&=\int_{B_{k}}\left(\int_{|y|\geq |x|}\frac{(b(x)-b(y))}{\omega(B(0,|y|))}f(y)\omega(y)dy\right)^p\omega(x)dx\\
&\leq \int_{B_{k}}\left(\int_{2^{k}\geq |y|\geq |x|}\frac{|b(x)-b(y)|}{\omega(B(0,|y|))}|f(y)|\omega(y)dy\right)^p\omega(x)dx\\
&+\int_{B_{k}}\left(\int_{|y|>2^{k}}\frac{|b(x)-b(y)|}{\omega(B(0,|y|))}|f(y)|\omega(y)dy\right)^p\omega(x)dx\\
&=:J_1+J_2.
\end{align*}

Similar to the estimates of $\mathcal{H}_{\omega,b}$, we have

\begin{align*}
J_1&\leq \int_{B_{k}}\left(\int_{ |y|\leq 2^{kn}}\frac{|b(x)-b(y)|}{\omega(B(0,|y|))}|f(y)|\omega(y)dy\right)^p\omega(x)dx\\
&\leq \int_{B_{k}}\left(\sum\limits_{i=-\infty}^{k}\int_{ B_i}\frac{|b(x)-b(y)|}{\omega(B(0,|y|))}|f(y)|\omega(y)dy\right)^p\omega(x)dx\\
&\leq  \int_{B_{k}}\left(\sum\limits_{i=-\infty}^{k}\frac{1}{\omega{(B_i)}}\int_{ B_i}|b(x)-b(y)||f(y)|\omega(y)dy\right)^p\omega(x)dx\\
&\leq \|b\|^p_{\dot{\mathcal{C}}^{p_1,\lambda_1}_\omega(\Bbb{R}^n)}\|f\|^p_{\dot{\mathcal{M}}^{p_2,\lambda_2}_\omega(\Bbb{R}^n)}\omega(B_{k})^{1+\lambda p}.
\end{align*}

Thus, it remains to give the estimates of $J_2$. By Lemma 2.4, we can decompose $J_2$ as
\begin{align*}
J_2
&\leq C \int_{B_{k}}\left(\sum\limits_{i=k}^\infty\frac{1}{\omega(B_i)}\int_{C_i}|b(x)-b(y)||f(y)|\omega(y)dy\right)^p\omega(x)dx\\
&\leq C\int_{B_{k}}\left(\sum\limits_{i=k}^\infty\frac{1}{\omega(B_i)}\int_{B_i}|b(x)-b_{B_{k},\omega}||f(y)|\omega(y)dy\right)^p\omega(x)dx\\
&+C \int_{B_{k}}\left(\sum\limits_{i=k}^\infty\frac{1}{\omega(B_i)}\int_{B_i}|b(y)-b_{B_i,\omega}||f(y)|\omega(y)dy\right)^p\omega(x)dx\\
&+C\int_{B_{k}}\left(\sum\limits_{i=k}^\infty\frac{1}{\omega(B_i)}\int_{B_i}|b_{B_i,\omega}-b_{B_{k},\omega}||f(y)|\omega(y)dy\right)^p\omega(x)dx\\
&=:CJ_{21}+CJ_{22}+CJ_{23}.
\end{align*}

For $J_{21}$,  we show the following fact.

\begin{equation}\label{resolution}
\sum\limits_{i=k}^\infty\frac{1}{\omega(B_i)}\int_{B_i}|f(y)|\omega(y)dy
\leq C\|f\|_{\dot{\mathcal{M}}^{p_2,\lambda_2}_\omega(\Bbb{R}^n)}\left(\omega(B_k)\right)^{\lambda_2}
\end{equation}
with   $y\in B_k$, $1<p_2<\infty$ and $\lambda_2<0$.

Then , using Lemma 2.3 and the H\"{o}lder inequality,   we get

\begin{align*}
&\sum\limits_{i=k}^\infty\frac{1}{\omega(B_i)}\int_{B_i}|f(y)|\omega(y)dy\\
&\leq C\sum\limits_{i=k}^\infty\frac{1}{\omega(B_i)}\left(\int_{B_i}|f(y)|^{p_2}\omega(y)dy\right)^{1/p_2}\left(\int_{B_i}\omega(y)^{(1-\frac{1}{p_2})p_2'}dy\right)^{1/p_2'}\\
&\leq C\|f\|_{\dot{\mathcal{M}}^{p_2,\lambda_2}_\omega(\Bbb{R}^n)}\sum\limits_{i=k}^\infty\omega(B_i)^{\lambda_2}\\
&\leq C\|f\|_{\dot{\mathcal{M}}^{p_2,\lambda_2}_\omega(\Bbb{R}^n)}\sum\limits_{i=k}^\infty\left({D_2}^{(i-k)}\omega(B_k)\right)^{\lambda_2} \\
&\leq C\|f\|_{\dot{\mathcal{M}}^{p_2,\lambda_2}_\omega(\Bbb{R}^n)}\left(\omega(B_k)\right)^{\lambda_2} \\
\end{align*}
and we finish the proof of (7).

Using Lemma 2.3,  the H\"{o}lder inequality and the conditions of Theorem 1.1,  we have

\begin{align*}
J_{21}
&\leq \int_{B_{k}}|b(x)-b_{B_{k,\omega}}|^p\omega(x)dx\left(\sum\limits_{i=k}^\infty\frac{1}{\omega(B_i)}\int_{B_i}\omega(y)|f(y)|dy\right)^p\\
&\leq C\left(\int_{B_{k}}|b(x)-b_{B_{k,\omega}}|^{p_1}\omega(x)dx\right)^{p/p_1}\left(\int_{B_{k}}\omega(x)^{(1-\frac{p}{p_1})(\frac{p_1}{p})'}dx\right)^{\frac{1}{(\frac{p_1}{p})'}}\\
&\times \left(\|f\|_{\dot{\mathcal{M}}^{p_2,\lambda_2}_\omega(\Bbb{R}^n)}\left(\omega(B_k)\right)^{\lambda_2}\right)^p\\
&\leq C\|b\|^p_{\dot{\mathcal{C}}^{p_1,\lambda_1}_\omega(\Bbb{R}^n)}\|f\|^p_{\dot{\mathcal{M}}^{p_2,\lambda_2}_\omega(\Bbb{R}^n)}\omega(B_{k})^{1+\lambda_1p}\left(\omega(B_k)^{\lambda_2}\right)^p\\
&\leq  C\|b\|^p_{\dot{\mathcal{C}}^{p_1,\lambda_1}_\omega(\Bbb{R}^n)}\|f\|^p_{\dot{\mathcal{M}}^{p_2,\lambda_2}_\omega(\Bbb{R}^n)}\omega(B_{k})^{1+\lambda p}.
\end{align*}

For $J_{22}$,  we need to show the following inequality.
\begin{equation}\label{resolution}
\sum\limits_{i=k}^{\infty}\frac{1}{\omega(B_{i})}\int_{B_i}|b(y)-b_{B_i,\omega}||f(y)|\omega(y)dy
\leq C\|b\|_{\dot{\mathcal{C}}^{p_1,\lambda_1}_\omega(\Bbb{R}^n)}\|f\|_{\dot{\mathcal{M}}^{p_2,\lambda_2}_\omega(\Bbb{R}^n)}\left(\omega(B_k)\right)^{\lambda}. \end{equation}
with $x\in B_k$,  $1<p_i<\infty$, $-\frac{1}{p_i}<\lambda_i<0$  and  $\lambda=\sum\limits_{i=1}^2\lambda_i$.

Denote $s$ by $1/s=1-1/{p_1}-1/{p_2}$. Then,  using  Lemma 2.3 and the H\"{o}lder inequality, we get

\begin{align*}
&\sum\limits_{i=k}^{\infty}\frac{1}{\omega(B_{i})}\int_{B_i}|b(y)-b_{B_i,\omega}||f(y)|\omega(y)dy\\
&\leq C\sum\limits_{i=k}^{\infty}\frac{1}{\omega(B_{i})}\left(\int_{B_i}|b(y)-b_{B_i,\omega}|^{p_1}\omega(y)dy\right)^{1/p_1}\left(\int_{B_i}|f(y)|^{p_2}\omega(y)dy\right)^{1/p_2}\left(\int_{B_i}\omega(y)^{\frac{1}{s}s}dy\right)^{1/s}\\
&\leq C\|b\|_{\dot{\mathcal{C}}^{p_1,\lambda_1}_\omega(\Bbb{R}^n)}\|f\|_{\dot{\mathcal{M}}^{p_2,\lambda_2}_\omega(\Bbb{R}^n)}\sum\limits_{i=k}^{\infty}\omega(B_i)^{\lambda}\\
&\leq C\|b\|_{\dot{\mathcal{C}}^{p_1,\lambda_1}_\omega(\Bbb{R}^n)}\|f\|_{\dot{\mathcal{M}}^{p_2,\lambda_2}_\omega(\Bbb{R}^n)} \sum\limits_{i=k}^{\infty}\left({D_2}^{(i-k)}\omega(B_k)\right)^{\lambda}\\
&\leq C\|b\|_{\dot{\mathcal{C}}^{p_1,\lambda_1}_\omega(\Bbb{R}^n)}\|f\|_{\dot{\mathcal{M}}^{p_2,\lambda_2}_\omega(\Bbb{R}^n)}\omega(B_{k})^{\lambda}\\
\end{align*}
and we finish the proof of (8).

Thus, we obtain
\begin{align*}
J_{22}
&=\int_{B_{k}}\left(\sum\limits_{i=k}^\infty\frac{1}{\omega(B_i)}\int_{B_i}|b(y)-b_{B_i,\omega}||f(y)|\omega(y)dy\right)^p\omega(x)dx\\
&\leq C\int_{B_{k}}\left(\|b\|_{\dot{\mathcal{C}}^{p_1,\lambda_1}_\omega(\Bbb{R}^n)}\|f\|_{\dot{\mathcal{M}}^{p_2,\lambda_2}_\omega(\Bbb{R}^n)}\omega(B_{k})^{\lambda}\right)^p\omega(x)dx\\
& \leq C\|b\|^p_{\dot{\mathcal{C}}^{p_1,\lambda_1}_\omega(\Bbb{R}^n)}\|f\|^p_{\dot{\mathcal{M}}^{p_2,\lambda_2}_\omega(\Bbb{R}^n)}\omega(B_{k})^{1+\lambda p}.
\end{align*}

To estimate  $J_{23}$,  using the conditions of Theorem 1.1,  Lemmas 2.3 and 2.5,   we get
\begin{align*}
J_{23}
&\leq C\|b\|^p_{\dot{\mathcal{C}}^{p_1,\lambda_1}_\omega(\Bbb{R}^n)}\int_{B_{k}}\left(\sum\limits_{i=k}^\infty\omega(B_i)^{-1}\int_{B_i}\omega(B_{k})^{\lambda_1}f(y)\omega(y)dy\right)^p\omega(x)dx\\
&\leq C\|b\|^p_{\dot{\mathcal{C}}^{p_1,\lambda_1}_\omega(\Bbb{R}^n)}\int_{B_{k}}\left(\sum\limits_{i=k}^\infty\omega(B_i)^{-1}\omega(B_{k})^{\lambda_1}
 \left(\int_{B_i}|f(y)|^{p_2}\omega(y)dy\right)^{1/p_2}\right.\\
& \left. \ \ \ \times\left(\int_{B_i}\omega(y)^{(1-\frac{1}{p_2})\frac{p_2}{p_2-1}}dy\right)^{\frac{p_2-1}{p_2}}\right)^p\omega(x)dx\\
&\leq C\|b\|^p_{\dot{\mathcal{C}}^{p_1,\lambda_1}_\omega(\Bbb{R}^n)}\int_{B_{k}}\left(\sum\limits_{i=k}^\infty\omega(B_i)^{-1}\omega(B_{k})^{\lambda_1}
 \|f\|_{\dot{\mathcal{M}}^{p_2,\lambda_2}_\omega(\Bbb{R}^n)}\omega(B_i)^{1+\lambda_2}\right)^p\omega(x)dx\\
& \leq
 C\|b\|^p_{\dot{\mathcal{C}}^{p_1,\lambda_1}_\omega(\Bbb{R}^n)}\|f\|^p_{\dot{\mathcal{M}}^{p_2,\lambda_2}_\omega(\Bbb{R}^n)}\int_{B_{k}}\omega(B_{k})^{\lambda_1p}
 \left(\sum\limits_{i=k}^\infty\omega(B_i)^{\lambda_2}\right)^p\omega(x)dx\\
& \leq
 C\|b\|^p_{\dot{\mathcal{C}}^{p_1,\lambda_1}_\omega(\Bbb{R}^n)}\|f\|^p_{\dot{\mathcal{M}}^{p_2,\lambda_2}_\omega(\Bbb{R}^n)}\int_{B_{k}}\omega(B_{k})^{\lambda_1p}
 \left(\sum\limits_{i=k}^\infty {D_2}^{{(i-k)}{\lambda_2}}\omega(B_{k})^{\lambda_2}\right)^p\omega(x)dx\\
& \leq
 C\|b\|^p_{\dot{\mathcal{C}}^{p_1,\lambda_1}_\omega(\Bbb{R}^n)}\|f\|^p_{\dot{\mathcal{M}}^{p_2,\lambda_2}_\omega(\Bbb{R}^n)}\omega(B_{k})^{\lambda p}\int_{B_{k}}
 \left(\sum\limits_{i=k}^\infty {D_2}^{{(i-k)}{\lambda_2}}\right)^p\omega(x)dx\\
&\leq \|b\|^p_{\dot{\mathcal{C}}^{p_1,\lambda_1}_\omega(\Bbb{R}^n)}\|f\|^p_{\dot{\mathcal{M}}^{p_2,\lambda_2}_\omega(\Bbb{R}^n)}\omega(B_{k})^{1+\lambda p}.
\end{align*}

Combing the estimates of $J_1, J_{21}, J_{22}$ and $ J_{23}$, we get (4) and  finish the proof of Theorem 1.1.

\section{Proof of Theorem 1.2}

By Lemma 2.6,  it suffices to show that for any $k\in \Bbb{Z}$,  the following estimates hold.
\begin{equation}\label{resolution}\int_{B_{k}}|\mathcal{H}_{\omega,b}(f)(x)|^p\omega(x)dx\leq C\omega(B_{k})^{1+2\lambda p}\|b\|^p_{\dot{\mathcal{C}}^{p,\lambda}_\omega(\Bbb{R}^n)}\|f\|^p_{\dot{\mathcal{M}}^{p,\lambda}_\omega(\Bbb{R}^n)},\end{equation}

\begin{equation}\label{resolution}\int_{B_{k}}|\mathcal{H}^\ast_{\omega,b}(f)(x)|^p\omega(x)dx\leq C\omega(B_{k})^{1+2\lambda p}\|b\|^p_{\dot{\mathcal{C}}^{p,\lambda}_\omega(\Bbb{R}^n)}\|f\|^p_{\dot{\mathcal{M}}^{p,\lambda}_\omega(\Bbb{R}^n)}.\end{equation}

To prove  (9),  using  Lemma 2.4,  we get
\begin{align*}
\int_{B_{k}}&|\mathcal{H}_{\omega,b}f(x)|^p\omega(x)dx
\leq \int_{B_{k}}\left(\frac{1}{\omega(B(0,|x|))}\int_{|y|<|x|}|b(x)-b(y)||f(y)|\omega(y)dy\right)^p\omega(x)dx\\
&\leq C\sum\limits_{j=-\infty}^{k}\int_{C_j}\left(\frac{1}{\omega(B(0,|x|))}\sum\limits_{i=-\infty}^j\int_{B_i}|b(x)-b_{B_j,\omega}||f(y)|\omega(y)dy\right)^p\omega(x)dx\\
&+C\sum\limits_{j=-\infty}^{k}\int_{C_j}\left(\frac{1}{\omega(B(0,|x|))}\sum\limits_{i=-\infty}^j\int_{B_i}|b(y)-b_{B_j,\omega}||f(y)|\omega(y)dy\right)^p\omega(x)dx\\
&\leq C\sum\limits_{j=-\infty}^{k}\int_{C_j}\frac{1}{\omega(B_j)^p}\left(\sum\limits_{i=-\infty}^j\int_{B_i}|b(x)-b_{B_j,\omega}||f(y)|\omega(y)dy\right)^p\omega(x)dx\\
&+C\sum\limits_{j=-\infty}^{k}\int_{C_j}\frac{1}{\omega(B_j)^p}\left(\sum\limits_{i=-\infty}^j\int_{B_i}|b(y)-b_{B_j,\omega}||f(y)|\omega(y)dy\right)^p\omega(x)dx\\
&=:CL_1+CL_2.
\end{align*}

For $L_1$, recall that $2<p<\infty$ and  $-\frac{1}{2p}<\lambda<0$. Moreover, in this case,  $p_2=p$,  $\lambda_2=\lambda$.  Then, using (5) in Section3, we have  $1<p<\infty$ and   $-1<\lambda<0$. Thus, we obtain
$$
\sum\limits_{i=-\infty}^j\int_{B_i}|f(y)|\omega(y)dy
\leq C\|f\|_{\dot{\mathcal{M}}^{p,\lambda}_\omega(\Bbb{R}^n)}\left(\omega(B_j)\right)^{\lambda+1},
$$

From  Lemma 2.3 and the H\"{o}lder inequality, there is

\begin{align*}
L_1
&\leq C\sum\limits_{j=-\infty}^{k}\frac{1}{\omega(B_j)^p}\int_{B_j}|b(x)-b_{B_j,\omega}|^p\omega(x)dx\left(\sum\limits_{i=-\infty}^j\int_{B_i}|f(y)|\omega(y)dy\right)^p\\
&\leq C\sum\limits_{j=-\infty}^{k}\frac{1}{\omega(B_j)^p}\left(\omega(B_j)^{\left(\lambda+\frac{1}{p}\right)}
  \|b\|_{\dot{\mathcal{C}}^{p,\lambda}_\omega(\Bbb{R}^n)}\right)^p
\left(\|f\|_{\dot{\mathcal{M}}^{p,\lambda}_\omega(\Bbb{R}^n)}\left(\omega(B_j)\right)^{\lambda+1}\right)^ p\\
&\leq C\|b\|^p_{\dot{\mathcal{C}}^{p,\lambda}_\omega(\Bbb{R}^n)}\|f\|^p_{\dot{\mathcal{C}}^{p,\lambda}_\omega(\Bbb{R}^n)}
  \sum\limits_{j=-\infty}^{k}\left(\omega(B_j)\right)^{\left(1+2\lambda p\right)} \\
&\leq C\|b\|^p_{\dot{\mathcal{C}}^{p,\lambda}_\omega(\Bbb{R}^n)}\|f\|^p_{\dot{\mathcal{C}}^{p,\lambda}_\omega(\Bbb{R}^n)}
  \sum\limits_{j=-\infty}^{k}\left({D_2}^{\left(j-k\right)}\omega(B_{k})\right)^{\left(1+2\lambda p\right)}\\
&\leq C\|b\|^p_{\dot{\mathcal{C}}^{p,\lambda}_\omega(\Bbb{R}^n)}\|f\|^p_{\dot{\mathcal{C}}^{p,\lambda}_\omega(\Bbb{R}^n)}
  \left(\omega(B_{k})\right)^{\left(1+2\lambda p\right)}\\
\end{align*}

For $L_2$,  using  (6) in Section 3 and it is easy to see

$$
\sum\limits_{i=-\infty}^j\int_{B_i}|b(y)-b_{B_j,\omega}||f(y)|\omega(y)dy
\leq C\|b\|_{\dot{\mathcal{C}}^{p,\lambda}_\omega(\Bbb{R}^n)}\|f\|_{\dot{\mathcal{M}}^{p,\lambda}_\omega(\Bbb{R}^n)}\left(\omega(B_j)\right)^{1+2\lambda}.
$$

Then, using Lemma 2.3, the H\"{o}lder inequality and the conditions of Theorem 1.2, we have

\begin{align*}
L_2
&=\sum\limits_{j=-\infty}^{k}\int_{C_j}\frac{1}{\omega(B_j)^p}\left(\sum\limits_{i=-\infty}^j\int_{B_i}|b(y)-b_{B_j,\omega}||f(y)|\omega(y)dy\right)^p\omega(x)dx\\
&\leq C\|b\|^p_{\dot{\mathcal{C}}^{p,\lambda}_\omega(\Bbb{R}^n)}\|f\|^p_{\dot{\mathcal{C}}^{p,\lambda}_\omega(\Bbb{R}^n)}\sum\limits_{j=-\infty}^{k}\int_{C_j}\frac{1}{\omega(B_j)^p}\left(\omega(B_j)^{1+2\lambda}\right)^p\omega(x)dx\\
&\leq C \|b\|^p_{\dot{\mathcal{C}}^{p,\lambda}_\omega(\Bbb{R}^n)}\|f\|^p_{\dot{\mathcal{C}}^{p,\lambda}_\omega(\Bbb{R}^n)}\sum\limits_{j=-\infty}^{k}\omega(B_j)^{1+2\lambda p}\\
&\leq C\|b\|^p_{\dot{\mathcal{C}}^{p,\lambda}_\omega(\Bbb{R}^n)}\|f\|^p_{\dot{\mathcal{C}}^{p,\lambda}_\omega(\Bbb{R}^n)}\omega(B_{k})^{1+2\lambda p}\sum\limits_{j=-\infty}^{k}D_2^{(j-k)(1+2\lambda p)}\\
&\leq C\|b\|^p_{\dot{\mathcal{C}}^{p,\lambda}_\omega(\Bbb{R}^n)}\|f\|^p_{\dot{\mathcal{C}}^{p,\lambda}_\omega(\Bbb{R}^n)}\omega(B_{k})^{1+2\lambda p}.
\end{align*}

Combing the estimates of $L_1$, $L_2$,   we find that (9) is true.

Now,  let us  focus on the proof of (10).
 First,  we   decompose $\int_{B_{k}}|\mathcal{H}^\ast_{\omega,b}(f)(x)|^p\omega(x)dx$ as follows.

 \begin{align*}
\int_{B_{k}}|\mathcal{H}^\ast_{\omega,b}(f)(x)|^p\omega(x)dx&=\int_{B_{k}}\left(\int_{|y|\geq |x|}\frac{(b(x)-b(y))}{\omega(B(0,|y|))}f(y)\omega(y)dy\right)^p\omega(x)dx\\
&\leq \int_{B_{k}}\left(\int_{2^{k}\geq |y|\geq |x|}\frac{|b(x)-b(y)|}{\omega(B(0,|y|))}|f(y)|\omega(y)dy\right)^p\omega(x)dx\\
&+\int_{B_{k}}\left(\int_{|y|>2^{k}}\frac{|b(x)-b(y)|}{\omega(B(0,|y|))}|f(y)|\omega(y)dy\right)^p\omega(x)dx\\
&=:M_1+M_2.
\end{align*}

Similar to the estimates of $\mathcal{H}_{\omega,b}$, we have

\begin{align*}
M_1
&\leq \int_{B_{k}}\left(\int_{ |y|\leq 2^{k}}\frac{|b(x)-b(y)|}{\omega(B(0,|y|))}|f(y)|\omega(y)dy\right)^p\omega(x)dx\\
&\leq \int_{B_{k}}\left(\sum\limits_{i=-\infty}^{k}\int_{ C_i}\frac{|b(x)-b(y)|}{\omega(B(0,|y|))}|f(y)|\omega(y)dy\right)^p\omega(x)dx\\
&\leq C \int_{B_{k}}\left(\sum\limits_{i=-\infty}^{k}\frac{1}{\omega{(B_i)}}\int_{ C_i}|b(x)-b(y)||f(y)|\omega(y)dy\right)^p\omega(x)dx\\
&\leq C\|b\|^p_{\dot{\mathcal{C}}^{p,\lambda}_\omega(\Bbb{R}^n)}\|f\|^p_{\dot{\mathcal{M}}^{p,\lambda}_\omega(\Bbb{R}^n)}\omega(B_{k})^{1+2\lambda p}.
\end{align*}

Thus, it remains to give the estimates of $M_2$. Using Lemma 2.4,  we can decompose $M_2$ as
\begin{align*}
M_2
&\leq C \int_{B_{k}}\left(\sum\limits_{i=k}^\infty\frac{1}{\omega(B_i)}\int_{C_i}|b(x)-b(y)||f(y)|\omega(y)dy\right)^p\omega(x)dx\\
&\leq C\int_{B_{k}}\left(\sum\limits_{i=k}^\infty\frac{1}{\omega(B_i)}\int_{B_i}|b(x)-b_{B_{k},\omega}||f(y)|\omega(y)dy\right)^p\omega(x)dx\\
& \ \ \ + C\int_{B_{k}}\left(\sum\limits_{i=k}^\infty\frac{1}{\omega(B_i)}\int_{B_i}|b(y)-b_{B_{i},\omega}||f(y)|\omega(y)dy\right)^p\omega(x)dx\\
& \ \ \
+C\int_{B_{k}}\left(\sum\limits_{i=k}^\infty\frac{1}{\omega(B_i)}\int_{B_i}|b_{B_{i},\omega}-b_{B_{k},\omega}||f(y)|\omega(y)dy\right)^p\omega(x)dx\\
&=:CM_{21}+CM_{22}+CM_{23}.
\end{align*}

For $M_{21}$,  using (7) in Section 3 and the conditions of Theorem 1.2, we may easily get
$$
\sum\limits_{i=k}^\infty\frac{1}{\omega(B_i)}\int_{B_i}|f(y)|\omega(y)dy
\leq\|f\|_{\dot{\mathcal{M}}^{p,\lambda}_\omega(\Bbb{R}^n)}\left(\omega(B_k)\right)^{\lambda}.
$$

Then, from Lemma 2.3 and the H\"{o}lder inequality,   we obtain
\begin{align*}
M_{21}
&\leq \int_{B_{k}}|b(x)-b_{B_{k,\omega}}|^p\omega(x)dx\left(\sum\limits_{i=k}^\infty\frac{1}{\omega(B_i)}\int_{B_i}|f(y)|\omega(y)dy\right)^p\\
&\leq \left(\omega(B_{k})^{\left(\lambda+\frac{1}{p}\right)}\|b\|_{\dot{\mathcal{C}}^{p,\lambda}_\omega(\Bbb{R}^n)}\right)^p \left(\|f\|_{\dot{\mathcal{M}}^{p,\lambda}_\omega(\Bbb{R}^n)}\left(\omega(B_k)\right)^{\lambda}\right)^p\\
&\leq  C\|b\|^p_{\dot{\mathcal{C}}^{p,\lambda}_\omega(\Bbb{R}^n)}\|f\|^p_{\dot{\mathcal{C}}^{p,\lambda}_\omega(\Bbb{R}^n)}\omega(B_{k})^{1+2\lambda p}.
\end{align*}

For  $M_{22}$,  using  (8) in Section 3 and we can easily obtain
$$
\sum\limits_{i=k}^{\infty}\frac{1}{\omega(B_{i})}\int_{B_i}|b(y)-b_{B_i,\omega}||f(y)|\omega(y)dy
\leq C\|b\|_{\dot{\mathcal{C}}^{p,\lambda}_\omega(\Bbb{R}^n)}\|f\|_{\dot{\mathcal{M}}^{p,\lambda}_\omega(\Bbb{R}^n)}\left(\omega(B_k)\right)^{2\lambda}.
$$

Then,  applying Lemma 2.3 and the H\"{o}lder inequality,  there is

\begin{align*}
M_{22}
&\leq C\int_{B_{k}}\left(\|b\|_{\dot{\mathcal{C}}^{p,\lambda}_\omega(\Bbb{R}^n)}\|f\|_{\dot{\mathcal{M}}^{p,\lambda}_\omega(\Bbb{R}^n)}\left(\omega(B_k)\right)^{2\lambda}\right)^p\omega(x)dx\\
&\leq C\|b\|^p_{\dot{\mathcal{C}}^{p,\lambda}_\omega(\Bbb{R}^n)}\|f\|^p_{\dot{\mathcal{C}}^{p,\lambda}_\omega(\Bbb{R}^n)}\omega(B_{k})^{1+2\lambda p}.
\end{align*}

Finally, we  give the estimates of   $M_{23}$.   using Lemma 2.5, the fact $i\geq k$ and $\lambda< 0$, there is
$$|b_{B_{k,\omega}}-b_{B_i,\omega}|\leq C\|b\|^p_{\dot{\mathcal{C}}^{p,\lambda}_\omega(\Bbb{R}^n)}\omega(B_{k})^{\lambda}.$$

Thus,  we have
\begin{align*}
&M_{23}
\leq C\|b\|^p_{\dot{\mathcal{C}}^{p,\lambda}_\omega(\Bbb{R}^n)}\omega(B_{k})^{\lambda p}\int_{B_{k}}\left(\sum\limits_{i=k}^\infty\frac{1}{\omega(B_i)}\int_{B_i}|f(y)|\omega(y)dy\right)^p\omega(x)dx\\
&\leq C\|b\|^p_{\dot{\mathcal{C}}^{p,\lambda}_\omega(\Bbb{R}^n)}\omega(B_{k})^{\lambda p}\\
&\times\int_{B_{k}}\left(\sum\limits_{i=k}^\infty\frac{1}{\omega(B_i)}\left(\int_{B_i}|f(y)|^{p}\omega(y)dy\right)^{1/p}\left(\int_{B_i}\omega(y)^{(1-\frac{1}{p})p'}dy\right)^{1/p'}\right)^p\omega(x)dx\\
&\leq C\|b\|^p_{\dot{\mathcal{C}}^{p,\lambda}_\omega(\Bbb{R}^n)}\|f\|^p_{\dot{\mathcal{C}}^{p,\lambda}_\omega(\Bbb{R}^n)}\omega(B_{k})^{\lambda p}\int_{B_{k}}\left(\sum\limits_{i=k}^\infty\omega(B_i)^{\lambda}\right)^p\omega(x)dx\\
&\leq C\|b\|^p_{\dot{\mathcal{C}}^{p,\lambda}_\omega(\Bbb{R}^n)}\|f\|^p_{\dot{\mathcal{C}}^{p,\lambda}_\omega(\Bbb{R}^n)}\omega(B_{k})^{{2\lambda } p}{\left(
\sum\limits_{i=k}^\infty{{D_2}^{(i-k){\lambda }}}\right)}^p\int_{B_{k}}\omega(x)dx\\
&\leq C\|b\|^p_{\dot{\mathcal{C}}^{p,\lambda}_\omega(\Bbb{R}^n)}\|f\|^p_{\dot{\mathcal{C}}^{p,\lambda}_\omega(\Bbb{R}^n)}\omega(B_{k})^{1+2\lambda p}.
\end{align*}

Combing the estimates of $M_1, M_{21}, M_{22}, M_{23}$,   we complete the proof of (10) and the proof of Theorem 1.2 has been finished.

\section{Proof of Theorem 1.3.}

We just give the proof of $(a_3)\Rightarrow (b_3)$. From (2) and the H\"{o}lder inequality, we have
\begin{align*}
&\frac{1}{\omega(B)^{1+p_1\lambda_1}}\int_B|b(y)-b_{B,\omega}|^{p_1}\omega(y)dy\leq \frac{C}{\omega(B)^{p_1\lambda_1}}\sup\limits_{y\in B}|b(y)-b_{B,\omega}|^{p_1}\\
&\leq \frac{C}{\omega(B)^{p_1\lambda_1}}\left(\frac{1}{\omega(B)}\int_B|b(y)-b_{B,\omega}|\omega(y)dy\right)^{p_1}\\
&\leq \frac{C}{\omega(B)^{p_1\lambda_1}}\left[\frac{1}{\omega(B)}\left(\int_B|b(y)-b_{B,\omega}|^p\omega(y)dy\right)^{1/p}\left(\int_B\omega(y)^{p'(1-\frac{1}{p})}dy\right)^{1/p'}\right]^{p_1}\\
&\leq \frac{C}{\omega(B)^{p_1\lambda_1}}\left[\frac{1}{\omega(B)}\left(\int_B|b(y)-b_{B,\omega}|^p\omega(y)dy\omega(B)^{1/p'}\right)^{1/p}\right]^{p_1}\\
&=\frac{C}{\omega(B)^{p_1\lambda_1}}\left(\frac{1}{\omega(B)}\int_B|b(y)-b_{B,\omega}|^p\omega(y)dy\right)^{\frac{p_1}{p}}.
\end{align*}
As
\begin{align*}
&\int_B|b(y)-b_{B,\omega}|^p\omega(y)dy\\
&=\int_B\left|b(y)-\frac{1}{\omega(B)}\int_Bb(z)\omega(z)dz\right|^p\omega(y)dy\\
&=\int_B\left|\frac{1}{\omega(y)}\int_B(b(y)-b(z))\omega(z)dz\right|^p\omega(y)dy\\
&\leq \frac{1}{\omega(B)^p}\int_B\left|\int_B|b(y)-b(z)|\omega(z)dz\right|^p\omega(y)dy\\
&\leq \frac{1}{\omega(B)^p}\int_B\left|\int_{\{z\in B:|z|<|y|\}}\chi_B(z)|b(y)-b(z)|\omega(z)dz\right|^p\omega(y)dy\\
&+ \frac{1}{\omega(B)^p}\int_B\left|\int_{\{z\in B:|z|\geq |y|\}}\chi_B(z)|b(y)-b(z)|\omega(z)dz\right|^p\omega(y)dy\\
&=:A+B.
\end{align*}
For $A$, as $y\in B$ implies $\omega(B(0,|y|))\leq \omega(B(0,r))=:\omega(B)$, there is
\begin{align*}
A&=\frac{1}{\omega(B)^p}\int_B\omega(B(0,|y|))^p\left|\frac{1}{\omega(B(0,|y|))}\int_{\{z\in B:|z|<|y|\}}\chi_B(z)|b(y)-b(z)|\omega(z)dz\right|^p\omega(y)dy\\
&=\frac{1}{\omega(B)^p}\int_B\omega(B(0,|y|))^p|\mathcal{H}_{\omega,|b|}(\chi_B)(y)|^p\omega(y)dy\\
&\leq \int_B|\mathcal{H}_{\omega,|b|}(\chi_B)(y)|^p\omega(y)dy\\
&\leq \omega(B)^{1+\lambda p}\|\mathcal{H}_{\omega,|b|}(\chi_B)\|^p_{\dot{\mathcal{M}}^{p,\lambda}_\omega}\\
&\leq \omega(B)^{1+\lambda p}\|\chi_B\|^p_{\dot{\mathcal{M}}^{p_2,\lambda_2}_\omega}\leq \omega(B)^{1+\lambda_1p}.
\end{align*}
For $B$, as $z\in B$ implies $\omega(B(0,|z|))\leq \omega(B(0,r))=:\omega(B)$, we obtain
\begin{align*}
B&=\frac{1}{\omega(B)^p}\int_{B}\left|\int_{z\in B:|z|\geq |y|}\frac{\chi_B(z)|b(y)-b(z)|\omega(z)}{\omega(B(0,|z|))}\omega(B(0,|z|))dz\right|^p\omega(y)dy\\
&\leq \int_B|\mathcal{H}^\ast_{\omega,|b|}(\chi_B)(y)|^p\omega(y)dy\\
&\leq \omega(B)^{1+\lambda p}\|\mathcal{H}^\ast_{\omega,|b|}(\chi_B)\|^p_{\dot{\mathcal{M}}^{p,\lambda}_\omega}\\
&\leq \omega(B)^{1+\lambda p}\|\chi_B\|^p_{\dot{\mathcal{M}}^{p_2,\lambda_2}_\omega}\leq \omega(B)^{1+\lambda_1p}.
\end{align*}
Combing the estimates of $A$ and $B$, it is easy to see
$$\frac{1}{\omega(B)^{1+p_1\lambda_1}}\int_B|b(y)-b_{B,\omega}|^{p_1}\omega(y)dy\leq C,$$
which implies $b\in \dot{\mathcal{C}}^{p_1,\lambda_1}_\omega(\Bbb{R}^n)$.

\section{Proof of Theorem 1.4.}

We just give the proof of $(a_4)\Rightarrow (b_4)$.  For any ball $B=B(0,R)$, we have
\begin{align*}
&\frac{1}{\omega(B)^{1+\lambda p}}\int_B|b(y)-b_{B,\omega}|^p\omega(y)dy\\
&=\frac{1}{\omega(B)^{1+\lambda p}}\int_B\left|\frac{1}{\omega(B)}\int_B(b(y)-b(z))\omega(z)dz\right|^p\omega(y)dy\\
&\leq \frac{1}{\omega(B)^{1+\lambda p+p}}\int_B\left|\int_{\{z\in B:|z|< |y|\}}|b(y)-b(z)|\omega(z)dz\right|^p\omega(y)dy\\
&+ \frac{1}{\omega(B)^{1+\lambda p+p}}\int_B\left|\int_{\{z\in B:|z|\geq |y|\}}|b(y)-b(z)|\omega(z)dz\right|^p\omega(y)dy\\
&=:G+H.
\end{align*}

For $G$, by the boundedness of $\mathcal{H}_{\omega,|b|}$ from $\dot{\mathcal{M}}^{p,\lambda}_\omega(\Bbb{R}^n)$ to $\dot{\mathcal{M}}^{p,2\lambda}_\omega(\Bbb{R}^n)$, there is

\begin{align*}
G&\leq \frac{1}{\omega(B)^{1+\lambda p+p}}\int_B\omega(B(0,|y|))^p|\mathcal{H}_{\omega,|b|}(\chi_B)(z)|^p\omega(y)dy\\
&\leq \frac{C}{\omega(B)^{1+\lambda p}}\int_{B}|\mathcal{H}_{\omega,|b|}(\chi_B)(z)|^p\omega(y)dy\\
&\leq  \frac{C}{\omega(B)^{1+\lambda p}}\|\mathcal{H}_{\omega,|b|}(\chi_B)(\cdot)\|_{\dot{\mathcal{M}}^{p,2\lambda}_\omega}^p\omega(B)^{p(1+2\lambda)}\\
&\leq \omega(B)^{p\lambda}\|\chi_B\|^p_{\dot{\mathcal{M}}^{p,\lambda}_\omega}\leq C.
\end{align*}
For $H$, using the fact  $\mathcal{H}^\ast_{\omega,|b|}$ is from $\dot{\mathcal{M}}^{p,\lambda}_\omega(\Bbb{R}^n)$ to $\dot{\mathcal{M}}^{p,2\lambda}_\omega(\Bbb{R}^n)$, we obtain
\begin{align*}
H&\leq \frac{1}{\omega(B)^{1+\lambda p+p}}\int_B\left|\int_{\{z\in B:|z|\geq |y|\}}\omega(B(0,|z|))\frac{|b(y)-b(z)|\omega(z)}{\omega(B(0,|z|))}dz\right|^p\omega(y)dy\\
&\leq \frac{1}{\omega(B)^{1+\lambda p}}\int_{B}|\mathcal{H}^\ast_{\omega,|b|}(\chi_B)(z)|^p\omega(y)dy\\
&\leq \frac{1}{\omega(B)^{1+\lambda p}}\|\mathcal{H}^\ast_{\omega,|b|}(\chi_B)(\cdot)\|_{\dot{\mathcal{M}}^{p,2\lambda}_\omega}^p\omega(B)^{p(1+2\lambda)}\\
&\leq \omega(B)^{p\lambda}\|\chi_B\|^p_{\dot{\mathcal{M}}^{p,\lambda}_\omega}\leq C.
\end{align*}
Consequently, the proof of Theorem 1.4 has been finished.

\vspace{1cm} \noindent {\bf {References}} \small

\begin{itemize}
\item [{[1]}] J. Alvarez, J. Lakey and M.
 Guzm\'{a}n-Partida, \textit{Spaces of bounded $\lambda$-central mean
oscilation, Morrey spaces, and $\lambda$-central Carleson measures},
Collect. Math.,  \textbf{51} (2000), 1-47.

\item [{[2]}] M. Christ, L. Grafakos,  \textit{ Best constants for two nonconvolution inequalities}, Proc. Amer. Math. Soc., \textbf{123} (1995),  1687-1693.

\item [{[3]}] R. Coifman and C. Fefferman, \textit{Weighted norm inequalities for maximal functions and singular integrals}, Studia Math., \textbf{51} (1974), 241-250.

\item [{[4]}]Z.W. Fu, Y. Lin and S.Z. Lu, \textit{$\lambda$-central BMO estimates for commutators of singular integrals with rough
kernels}, Acta Math. Sinica, English Ser., \textbf{24} (2008),  373-386.

\item [{[5]}] Z.W. Fu, Z. G. Liu, S.Z. Lu and H.B. Wang, \textit{Characterization for commutators of $n$-dimensional fractional Hardy operators}, Sci. China, Ser. A, \textbf{50} (2007), 418-426.

\item [{[6]}]L. Grafakos, \textit{Classical and Modern Fourier Analysis}, Pearson Education, Inc.,Upper Saddle River, New Jersey, 2004.

\item [{[7]}] G. Hardy,  \textit{Note on a theorem of Hilbert}, Math. Z., \textbf{6} (1920),   314-317.

\item [{[8]}] Y. Komori, \textit{Notes on singular integrals on some inhomogeneous herz
spaces}, Taiwan. J. Math., \textbf{8} (2004), 547-556.

\item [{[9]}] Y. Komori and S. Shirai, \textit{Weighted Morrey spaces and a singular integral operator},  Math. Nachr., \textbf{282} (2009), 219-231.

\item [{[10]}]S.Z. Lu and D.C. Yang, \textit{The Littlewood-Paley function and $\varphi$-transform characterization of a new Hardy space $HK_2$ associated with Herz space}. Studia Math., \textbf{101} (1992), 285-298.

\item [{[11]}] S.Z. Lu and D.C. Yang, \textit{The central BMO space and Littlewood operators}. Approx. Theory Appl., \textbf{11} (1995), 72-94.

\item [{[12]}]C.B. Morrey, On the solution of quasi-linear elliptic partial differential eqautions, Trans. Amer. Math. Soc., \textbf{43}(1938), 126-166.

\item [{[13]}]B. Muckenhoupt, \textit{Weighted norm inequalities for the Hardy maximal function}, Trans. Amer. Math. Soc., \textbf{165}(1972), 207-226.
\item [{[14]}] S.G. Shao, S.Z. Lu,\textit{Some characterizations of Campanato spaces via commutators on Morrey spaces}, Pacific J. Math., \textbf{264}(2013),  221-234.
\item [{[15]}]S.G. Shao, S.Z. Lu, \textit{A characterization of Campanato space via commutator of fractional integral},  J. Math. Anal. Appl., \textbf{419} (2014), no.1, 123-137.
\item [{[16]}]  S.G. Shao, S.Z. Lu, \textit{Characterization of the central Campanato space via the commutator operator of Hardy type}, J. Math. Anal. Appl., \textbf{429} (2015), 713-732.

\item [{[17]}] X. Yu and X.X. Tao, \textit{Boundedness for a class of generalized commutators on  $\lambda$-central Morrey space}, Acta Math. Sin. (Engl. Ser.), \textbf{29} (2013), 1917-1926.

\item [{[18]}]X. Yu, H.H. Zhang and G.P. Zhao, \textit{Weighted boundedness of some integral operators on weightd $\lambda$-central Morrey space}, Appl. Math., J. Chinese Univ., Ser. B, \textbf{31}(2016), 331-342.

\item [{[19]}] F.Y. Zhao, S.Z.Lu, \textit{A  characterization of $\lambda$-central BMO space}, Front. Math. China, \textbf{8} (2013), 229-238.

\end{itemize}\vskip 10mm

\end{document}